\documentclass[10pt,a4paper,twoside,twocolumn]{IEEEtran}

\usepackage{amsfonts}
\usepackage{amssymb}
\usepackage{amsmath}
\usepackage{graphicx}
\usepackage{subfigure}
\usepackage{psfrag}
\let\labelindent\relax\usepackage{enumitem}
\usepackage{algorithm}
\usepackage{cite}
\usepackage[left=1.9cm,top=2.0cm,right=1.9cm, bottom=2.0cm]{geometry}
\usepackage{balance}

\newcommand\abs[1]{\lvert#1\rvert}
\newcommand\diag[1]{\mathrm{diag}(#1)}

\begin{document}

\title{Dynamic Spectrum Sensing Through Accelerated Particle Swarm Optimization}

\author{Alexandros~E.~Paschos,~\IEEEmembership{Student Member,~IEEE},
Vasileios~M.~Kapinas,~\IEEEmembership{Member,~IEEE},
Georgia~D.~Ntouni,~\IEEEmembership{Student Member,~IEEE},
Leontios~J.~Hadjileontiadis,~\IEEEmembership{Senior Member,~IEEE},
and George~K.~Karagiannidis~\IEEEmembership{Fellow,~IEEE}
\thanks{A.~E.~Paschos is with Institut Sup\'{e}rieur de l'A\'{e}ronautique et de l'Espace, Universit\'{e} de Toulouse, 31400, France
(email: alexander.pasx@gmail.com).}%
\thanks{V. M. Kapinas, G. D. Ntouni, L. J. Hadjileontiadis, and G. K. Karagiannidis are with the Department of Electrical and Computer Engineering, Aristotle University of Thessaloniki, 54636, Greece
(email: \{kapinas, gntouni, leontios, geokarag\}@auth.gr).}%
\thanks{L.~J.~Hadjileontiadis is also with the Department of Electrical and Computer Engineering, Khalifa University of Science and Technology, Abu Dhabi, PO Box 127788, UAE
(email: leontios.h@kustar.ac.ae).}%
\thanks{The work of G. D. Ntouni has been supported by a scholarship from the Alexander Onassis Foundation.}%
}

\maketitle


\begin{abstract}
In this paper, a novel optimization algorithm, called the acceleration-aided particle swarm optimization (A-APSO), is proposed for reliable dynamic spectrum sensing in cognitive radio networks. In A-APSO, the acceleration variable of the particles in the swarm is also considered in the search space of the optimization problem. We show that the proposed A-APSO based spectrum sensing technique is more efficient in terms of performance than the corresponding one based on the standard particle swarm optimization algorithm.
\end{abstract}

\begin{IEEEkeywords}
Acceleration-aided particle swarm optimization, cognitive radio networks, cooperative spectrum sensing.
\end{IEEEkeywords}

\section{Introduction}\label{sec:intro}

Cognitive radio (CR) has emerged as a promising technology to give solutions to the continuously increasing traffic demand. Many techniques for efficient dynamic spectrum management in CR networks (CRNs) have been proposed so far. One of the most popular ones is the opportunistic use of spectrum, where unlicensed users are enabled to access licensed frequency bands detected to be idle \cite{Suraweera}. However, allowing secondary users to utilize licensed bands requires reliable spectrum sensing of weak primary signals~\cite{Stotas2,Nijsure,Stotas1}. To this end, cooperative spectrum sensing is a new design paradigm in CRNs that can provide significant multiplexing and diversity gains.

Particularly, cooperative transmission can greatly improve the spectrum access opportunity as well as sharing efficiency for cognitive users with the help of cooperative relay nodes \cite{Hanzo_2017}. Various algorithms have been proposed so far for the implementation of efficient cooperative spectrum sensing in CRNs, just to mention \cite{Quan, Chen_2014,Duong,daCosta,Zou,Lee,Manesh_2017,Alhammadi_2016,Dawoud1_2016,Dawoud2_2016,Zheng,Pradhan,Rashid,Udgata} and references therein. Among them, particle swarm optimization (PSO) has been recently proved to be a very handy technique for spectrum sensing and allocation in CRNs \cite{Manesh_2017,Alhammadi_2016,Dawoud1_2016,Dawoud2_2016,Zheng,Pradhan,Rashid,Udgata}. In a nutshell, the current research trend in the field is toward the determination of the optimum power allocation and the simplification of the relay selection process \cite{Hanzo_2017}.

Generally speaking, PSO is an effective computational method for optimizing continuous nonlinear functions \cite{Kennedy,Shi1,Shi2,Han_2017}. It is a simple, fast and efficient stochastic swarm intelligence algorithm used in many discrete optimization problems \cite{Pradhan}. PSO neither requires a differentiable objective function nor relies on a specific single variable initialization, while it is less complex than other evolutionary optimization methods, e.g., genetic algorithms \cite{Manesh_2017,Dawoud1_2016}. These merits render the PSO-based techniques attractive candidates for dealing with dynamic spectrum sensing in CRNs, which may involve non-convex and joint optimization of several parameters at the same time. Interestingly, it has been recently shown in \cite{Rashid} that a particularly tailored PSO algorithm is capable of further improving the computational complexity by considering the tradeoff between the detection performance and optimization time of the spectrum sensing process in CRNs.

The standard PSO method, originally proposed by Kennedy and Eberhart in \cite{Kennedy} and later refined by Shi and Eberhart in \cite{Shi1,Shi2}, has been applied in several scientific fields, such as in optimization analysis, computational intelligence, and scheduling applications. More than thirty PSO variants have been proposed so far to achieve accelerated results, just to mention \cite{Yang,Mohamed}. However, existing algorithms neglect the acceleration factor of the particles in the swarm, whereas they adopt the term ``accelerated" to characterize their convergence rate. Besides, in other techniques inspired from physics-based modeling, such as in \cite{Elkaim}, spring type forces among swarm particles may cause acceleration discontinuities or possible swarm splitting and overall anomalies in behavior.

In this paper, we propose a novel robust spectrum sensing technique for CRNs, employing an optimization algorithm, namely acceleration-aided PSO (A-APSO), inspired from recent applications in signal processing \cite{Apostolidis,Kaltsa}. Particularly, derived from physics laws, A-APSO enriches the swarm intelligence theory with the involvement of the acceleration factor of the swarm in the model equations. In the sequel, we describe in detail the A-APSO based spectrum sensing algorithm, and we finally compare its performance with the standard PSO-based method~\cite{Zheng}. Interestingly, in the proposed swarming model, the obstacles mentioned above are avoided due to the refined observation of the swarm. However, this advantage comes at the cost of some complexity increase in the algorithm, which though is affordable thanks to the continuously growing availability in computational resources.

\section{System Model}\label{sec:system}

\begin{figure}
\psfrag{H}[][][0.7]{\,$H_0$\,/\,$H_1$}%
\psfrag{C1}[][][0.8]{\,\,CR$_\text{1}$}%
\psfrag{C2}[][][0.8]{\,\,CR$_\text{2}$}%
\psfrag{CM}[][][0.8]{\,\,\,CR$_\text{M}$}%
\psfrag{L2}[][][0.75]{\,\,$\sum_{k=0}^{N-1}\abs{\cdot}^2$}
\psfrag{-//-}[][][0.75]{\,\,$\sum_{k=0}^{N-1}\abs{\cdot}^2$}%
\psfrag{S}[][][2.0]{$\Sigma$}%
\psfrag{D}[][][0.8]{\,Decision}%
\psfrag{FC}[][t][0.8]{\quad\quad\quad Fusion Center}%
\psfrag{r1}[][][0.8]{$r_1(k)$}
\psfrag{r2}[][][0.8]{$r_2(k)$}
\psfrag{rm}[][][0.7]{\,\,$r_M(k)$}
\psfrag{u1}[][][0.9]{\,$u_1$}
\psfrag{u2}[][][0.9]{\,$u_2$}\emph{}
\psfrag{um}[][][0.85]{\,\,\,$u_M$}
\psfrag{z1}[][][0.9]{\,$z_1$}
\psfrag{z2}[][][0.9]{\,$z_2$}
\psfrag{zm}[][][0.9]{\,\,\,\,$z_M$}
\psfrag{y1}[][][0.9]{$y_1$}
\psfrag{y2}[][][0.9]{$y_2$}
\psfrag{ym}[][][0.9]{\,\,$y_M$}
\psfrag{w1}[][][0.9]{\,\,$w_1$}
\psfrag{w2}[][][0.9]{\,\,$w_2$}
\psfrag{wm}[][][0.9]{\,\,\,\,\,$w_M$}
\psfrag{yfc}[][][0.9]{\,\,$y_{fc}$}
\centering%
{\includegraphics[width=1.0\columnwidth,trim=0 0 0 0,clip=true]{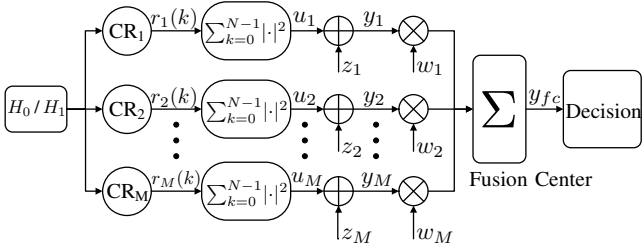}}
\caption{Weighting cooperation for dynamic spectrum sensing in a CRN.}\label{fig:crn}
\end{figure}

We consider the cooperative spectrum sensing network of Fig.~\ref{fig:crn}, consisting of $M$ CRs that send their locally sensed statistics to a fusion center. The binary hypothesis test, with $H_0$ and $H_1$ representing the hypothesis of signal being absent and present, respectively, at the $k$th time instant, is given~by
\begin{align}
&H_0:\,{r_{m}(k)}={n_{m}(k)},\quad m=1,2,\dotsc,M,\label{eq:hypotest0}\\
&H_1:\,{r_{m}(k)}={h_{m}s(k)}+n_{m}(k),\quad m=1,2,\dotsc,M,\label{eq:hypotest1}
\end{align}
where $s(k)$ is the transmitted signal from the primary user, $h_m$ and $n_m(k)\sim\mathcal{CN}(0,\sigma_m^2)$ denote the channel gain (being constant during the detection interval) and the zero-mean additive white Gaussian noise (AWGN) of variance $\sigma_m^2$, respectively, while $r_m(k)$ represents the received signal by the the $m$th CR.

The $m$th CR (or $m$th secondary user) calculates the statistic
\begin{align}\label{eq:userstat}
u_m=\sum_{k=0}^{N-1}{\abs{r_m(k)}}^2
\end{align}
considering a detection interval of $N$ samples. The statistics received by the fusion center through the control channel are
\begin{align}\label{eq:fusionstat}
y_m=u_m + z_m,\quad m=1,2,\dotsc,M,
\end{align}
where $z_m\sim\mathcal{CN}(0,\delta_m^2)$ is zero-mean AWGN of variance $\delta_m^2$ \cite{Quan}. The fusion center assigns non-negative weights $w_m\leq 1$ for $m=1,2,\dotsc,M$ and calculates the global statistic
\begin{align}\label{eq:globalstat}
y_{fc}=\sum_{m=1}^M{w_{m}y_{m}}=\mathbf{w}^T\mathbf{y},
\end{align}
where $\mathbf{w}={[w_1\cdots w_M ]}^T$ is the weight vector applied, and $\mathbf{y}={[y_1\cdots y_M ]}^T$ is the received vector by the fusion center.

\section{Problem Formulation}\label{sec:problem}

Given a targeted probability of false alarm $P_f$, the probability of detection $P_d$ can be calculated in terms of the transmitted power, $E_s=\sum_{k=0}^{N-1}{\abs{s(k)}}^2$ (assumed to be known at the fusion center of the CRN), the weight vector $\mathbf{w}$, and the channel power gain vector $\mathbf{g}=\left[\abs{h_1}^2\cdots\abs{h_M}^2\right]^T$~as~\cite{Zou}
\begin{align}\label{eq:detectprob}
 P_d=Q\left[\frac{Q^{-1}(P_f)\sqrt{{\mathbf{w}}^T\mathbf{a}\mathbf{w}}-{E_s{\mathbf{g}}^T{\mathbf{w}}}} {\sqrt{{\mathbf{w}}^T{\mathbf{b}\mathbf{w}}}}\right],
\end{align}
where $Q(x)=(1/\sqrt{2\pi})\int_{x}^{\infty}\mathrm{e}^{-t^2/2}\mathrm{d}t$ is the Gaussian $Q$-function, and $\mathbf{a},\mathbf{b}$ are $M\times 1$ vectors given by
\begin{align}\label{eq:AandB}
\begin{aligned}
&\mathbf{a} = 2N\diag{\boldsymbol{\sigma}}^2 + \diag{\boldsymbol{\delta}},\\
&\mathbf{b} = 2N\diag{\boldsymbol{\sigma}}^2 + \diag{\boldsymbol{\delta}}+ 4E_s\diag{\mathbf{g}}\diag{\boldsymbol{\sigma}},
\end{aligned}
\end{align}
with $\boldsymbol{\sigma}=[\sigma_1^2\cdots\sigma_M^2]^T$ and $\boldsymbol{\delta}=[\delta_1^2\cdots\delta_M^2]^T$. In \eqref{eq:AandB}, $\diag{\boldsymbol{\chi}}$ denotes a square matrix with the entries of vector $\boldsymbol{\chi}$ being in its main diagonal and zeros elsewhere.

The goal is to find the optimal weight vector that maximizes the probability of detection in \eqref{eq:detectprob}, which is equivalent to the minimization of its argument, since $P_d(\cdot)$ is a monotonically decreasing function. Thus, the optimization problem becomes
\begin{align}
\mathbf{w_\mathrm{opt}}=\underset{\mathbf{w}}{\arg}\min f(\mathbf{w}),\label{eq:minprob}
\end{align}
where $f(\mathbf{w})$ is the fitness function to be optimized, given~by
\begin{align}
f(\mathbf{w})={}&\frac{Q^{-1}(P_f)\sqrt{{\mathbf{w}}^T\mathbf{a}\mathbf{w}}-E_s{\mathbf{g}}^T \mathbf{w}}{\sqrt{{\mathbf{w}}^T\mathbf{b}\mathbf{w}}}.\label{eq:funcarg}
\end{align}

\section{A-APSO based Spectrum Sensing Algorithm}\label{sec:apso}

In this section, we introduce novel swarm model equations involving also the acceleration factor of the particles. To begin with, from Newton's $2$nd law, and after adopting notation $f^t\triangleq f(t)$ for mapping $t\mapsto f(t)$, velocity and position can be written in terms of the acceleration factor $a^t$ in discrete form~as
\begin{align}
 &v^t = v^{t-1}+a^{t}t,\label{eq:velocity1}\\
 &x^t = x^{t-1}+v^{t}t+\frac{1}{2}a^{t}t^2.\label{eq:position1}
\end{align}

In standard PSO algorithm, let $\mathbf{v}_i^t={[v_{i,1}^t\cdots v_{i,D}^t]}^T$ and $\mathbf{x}_i^t={[x_{i,1}^t\cdots x_{i,D}^t]}^T$ be the velocity and position vectors of the $i$th particle at iteration $t$, where $D$ is the number of particle dimensions, and $i=1,\dotsc,S$, with $S$ being the size of the swarm. Let now $\mathbf{p}_i^t={[p_{i,1}^t\cdots p_{i,D}^t]}^T$ be the best (position vector) solution obtained from the $i$th particle up to iteration $t$, and $\mathbf{p}_b^t={[p_{b,1}^t\cdots p_{b,D}^t]}^T$ the best (position vector) solution obtained from $\mathbf{p}_i^t$ in the population at iteration $t$. In this case, we can adjust the velocity and position equations~as
\begin{align}
 &v_{i,d}^t = \omega v_{i,d}^{t-1}+c_1 \xi(p_{i,d}^{t-1}-x_{i,d}^{t-1})+c_2 \eta(p_{b,d}^{t-1}-x_{i,d}^{t-1}), \label{eq:velocity2}\\
 &x_{i,d}^t = x_{i,d}^{t-1}+v_{i,d}^t,\label{eq:position2}
\end{align}
where $c_1$, $c_2$ are acceleration coefficients, $\xi,\eta$ are random numbers uniformly distributed in the $[0,1]$ interval, $\omega$ is an inertia weight, and $d=1,\dotsc,D=M$ \cite{Shi1}. After comparison of the standard PSO equations and the updated ones, we notice the absence of the acceleration factor, which we define now~as
\begin{align}\label{eq:acceler}
 a_{i,d}^t\triangleq c_1\xi(p_{i,d}^{t-1}-x_{i,d}^{t-1})+c_2\eta(p_{b,d}^{t-1}-x_{i,d}^{t-1}),
\end{align}
where $\mathbf{a}_i^t={[a_{i,1}^t\cdots a_{i,D}^t]}^T$ is the acceleration vector of the $i$th particle at iteration $t$. It is worth mentioning that the two terms on the right side of \eqref{eq:acceler}, which represent the cognitive and social components \cite{Elkaim}, respectively, are similar to elastic forces $\mathbf{f}=k_c\mathbf{x}$, where $k_c$ is the elastic constant and $\mathbf{x}$ the distance vector from the center of mass (Hooke's law).

Therefore, the velocity and position equations, considering also $t=1$ in \eqref{eq:velocity1} and \eqref{eq:position1}, can now be expressed~as
\begin{align}
 &v_{i,d}^t = {\omega v_{i,d}^{t-1}+a_{i,d}^t},\label{eq:velocity}\\
 &x_{i,d}^t = {x_{i,d}^{t-1}+v_{i,d}^{t}+\frac{1}{2}a_{i,d}^t}.\label{eq:position}
\end{align}

The proposed A-APSO is described in detail in Algorithm~\ref{alg:labelA}, where the computation of the fitness values (FVs) can be performed with the aid of the procedure shown in Algorithm~\ref{alg:labelB}. Note that, the inverse $Q$-function $Q^{-1}(\cdot)$ is available as a built-in function in many mathematical software packages (e.g., in MATLAB it can be computed via the \verb"qfuncinv" function). Additionally, in step~\ref{alg:labelA:thing4} of Algorithm~\ref{alg:labelA}, we can introduce further constraints on the acceleration behavior of swarm particles.

\begin{algorithm}
\caption{Novel A-APSO based spectrum sensing}\label{alg:labelA}
\begin{enumerate}[label=(\roman*), ref=(\roman*), leftmargin=2em]
 \item \label{alg:labelA:thing1} Set $t=0$ and randomly generate $v_{i,d}^t\in[-v_\mathrm{max},v_\mathrm{max}]$ and $x_{i,d}^t\in[0,1]$, where $i=1,\dotsc,S$ and $v_\mathrm{max}$ is the maximum selected velocity.
 \item \label{alg:labelA:thing2} Compute FV for each particle in the swarm and set $\mathbf{p}_i^t={[x_{i,1}^t\cdots x_{i,D}^t]}^T$ and $\mathbf{p}_b^t={[x_{b,1}^t\cdots x_{b,D}^t]}^T$, where $b$ is the index of the particle with the highest (negative) value.
 \item \label{alg:labelA:thing3} Set $t=t+1$ and update velocity $v_{i,d}^t$ according to~\eqref{eq:velocity}. If $v_{i,d}^t>v_\mathrm{max}$, then $v_{i,d}^t=v_\mathrm{max}$, while, if $v_{i,d}^t<-v_\mathrm{max}$, then $v_{i,d}^t=-v_\mathrm{max}$.
 \item \label{alg:labelA:thing4} Define acceleration factor as in~\eqref{eq:acceler}.
 \item Update position $x_{i,d}^t$ according to~\eqref{eq:position}.
 \item Compute the FV for each particle in the population. For the $i$th particle, if its FV is greater than the FV of $\mathbf{p}_i^{t-1}$, then set $\mathbf{p}_i^t={[x_{i,1}^t\cdots x_{i,D}^t]}^T$, else $\mathbf{p}_i^t=\mathbf{p}_i^{t-1}$. If the $i$th particle's FV is greater than the FV of $\mathbf{p}_b^{t-1}$, then set $\mathbf{p}_b^t={[x_{i,1}^t\cdots x_{i,D}^t]}^T$. Else, if there is no particle with FV greater than the FV of $\mathbf{p}_b^{t-1}$, then set $\mathbf{p}_b^t=\mathbf{p}_b^{t-1}$.
 \item If $t$ equals the max iteration, terminate, else go to~\ref{alg:labelA:thing3}.
\end{enumerate}
\end{algorithm}

\begin{algorithm}
\caption{Computation of the FVs used in Algorithm~\ref{alg:labelA}}\label{alg:labelB}
\begin{tabular}{l}
 $\phantom{1}1:$\, Set $P_f$ and compute $Q^{-1}(P_f)$\\
 $\phantom{1}2:$\, Set $N$ and $E_s\leftarrow 0$\\
 $\phantom{1}3:$\, \textbf{for} $k=0:N-1$ \textbf{do}\\
 $\phantom{1}4:$\, \quad\: $s(k+1)=1$\\
 $\phantom{1}5:$\, \quad\: $E_{s}=E_{s}+s(k+1)^2$\\
 $\phantom{1}6:$\, \textbf{end for}\\
 $\phantom{1}7:$\, Set $\boldsymbol{\gamma}$ and calculate $\mathbf{g}=\boldsymbol{\gamma} N/{E_{s}}$\\
 $\phantom{1}8:$\, Set $\boldsymbol{\sigma},\boldsymbol{\delta}$ and calculate $\mathbf{a},\mathbf{b}$\\
 $\phantom{1}9:$\, \textbf{for} $i=1:S$ \textbf{do}\\
 $10:$\, \quad\: Set $\mathbf{w}$ (i.e., $\mathbf{p}_i^t$ from Algorithm~\ref{alg:labelA})\\
 \quad\quad\quad\: and computation of $f(\mathbf{w})$\\
 $11:$\, \textbf{end for}
\end{tabular}
\end{algorithm}

\section{Simulations and Discussions}\label{sec:sim}

Simulations were based on a targeted probability of false alarm equal to $P_f=0.1$. The rest of system model parameters used are the number of CRs $M=6$, the number of problem dimensions $D=6$, the noise variance vectors $\boldsymbol{\sigma}=\boldsymbol{\delta}=[1 \: 1 \: 1 \: 1 \: 1 \: 1]^T$, the number of summary statistics (defining also the length of the detection interval) $N=20$, and the transmitted primary signal $s(k)=1$ for each time instant $k$. The algorithm parameters used are the acceleration coefficients $c_1=c_2=2$, the inertia weight $\omega=1$, the velocity range $[-v_\mathrm{max},v_\mathrm{max}]=[-5,5]$ (with $v_\mathrm{max}$ being its maximum value), while $\xi,\eta$ are random numbers uniformly distributed within $[0,1]$. The received signal-to-noise ratios (SNRs) at the CRs are $\boldsymbol{\gamma}=[-2.7 \: -1.2 \: -4.4 \: -4.5 \: -6.7 \: -4.7]^T$ dBs. Finally, our studies revealed that tuning the swarm size $S$ on the problem at hand is of minor importance in finding optimal solutions, as shown in Table~\ref{table:swarmsize}. However, in our problem, the optimal value that has been selected taking into account both performance and optimization is~$S=30$.

\begin{table}
\centering%
\caption{Impact of Swarm Size on Optimal Solutions}\label{table:swarmsize}
\begin{tabular}{c|c}
\hline\hline
$S$	(number of particles) & $P_d$ (probability of detection)\\
\hline\hline
\phantom{0}10 & 0.9387\\
\hline
\phantom{0}30 & 0.9413\\
\hline
\phantom{0}40 & 0.9413\\
\hline
\phantom{0}50 & 0.9413\\
\hline
\phantom{0}70 & 0.9413\\
\hline
100 & 0.9411\\
\hline\hline
\end{tabular}
\end{table}

\begin{figure}
\psfrag{Number of iterations}[t][][0.8][0]{Number of iterations}%
\psfrag{Probability of detection}[b][][0.8][0]{$P_d$}%
\centering%
\subfigure[Probability of detection for the proposed A-APSO ($\varepsilon=0.5$).\label{fig:apso}]%
{\includegraphics[width=1.0\columnwidth,trim=0 0 0 0,clip=true]{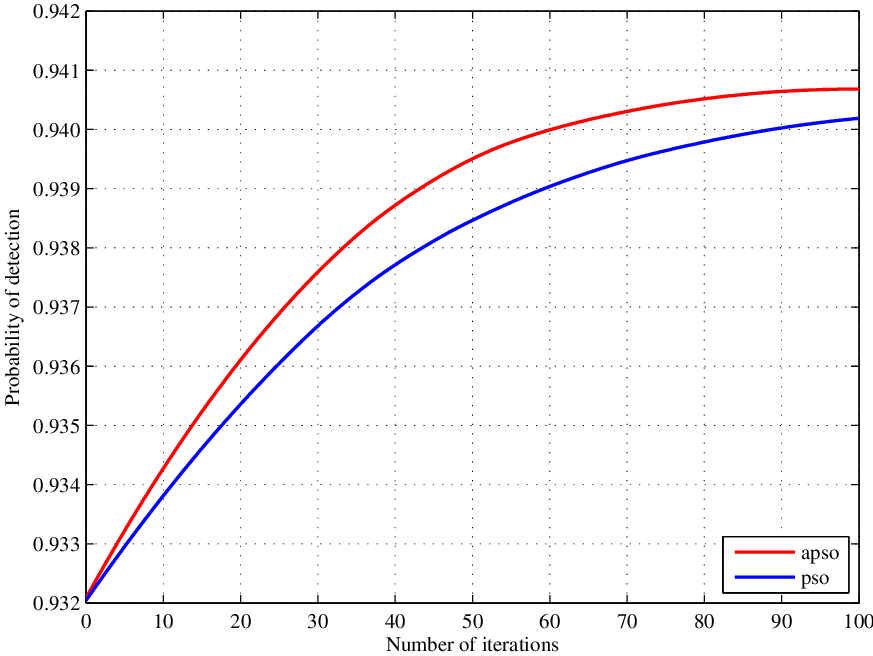}}
\subfigure[Probability of detection for the proposed A-APSO ($\varepsilon=2$).\label{fig:apsofurther}]%
{\includegraphics[width=1.0\columnwidth,trim=0 0 0 0,clip=true]{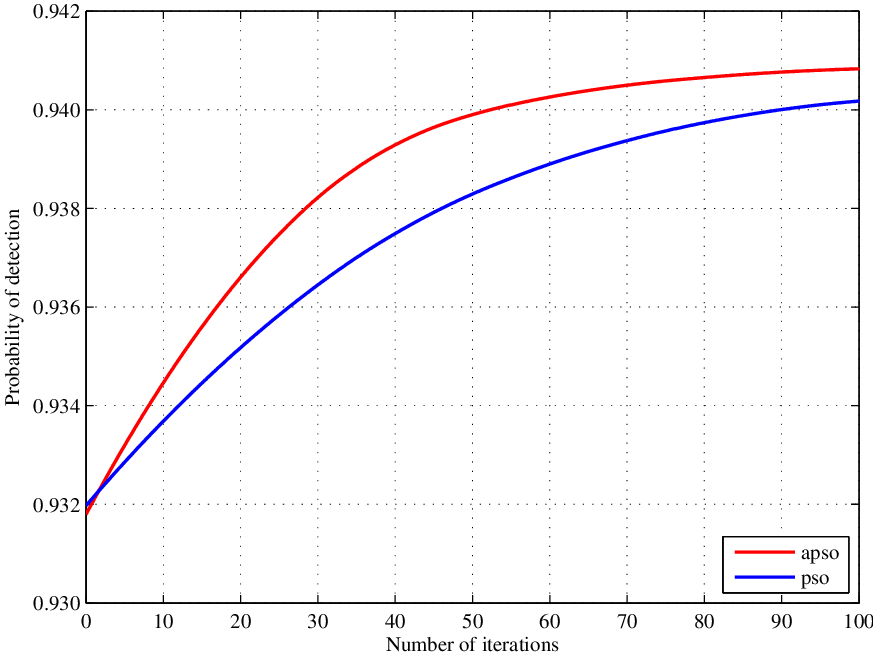}}
\caption{Performance comparison between PSO and A-APSO in terms of the probability of detection, $P_d$, vs. the number of iterations ($1000$ realizations per iteration and linear interpolation with $100$ query points).}\label{fig:apsodetection}
\end{figure}

The performance of the proposed A-APSO over the standard PSO algorithm in terms of the optimization problem described in Section~\ref{sec:problem} is shown in Fig.~\ref{fig:apsodetection}. The curves have been extracted after $1000$ realizations per iteration and linear interpolation with $100$ query points. From Fig.~\ref{fig:apso}, we can see that the proposed A-APSO algorithm has improved performance compared to the standard PSO over a number of $100$ iterations. Further studies concerning the equation of position in \eqref{eq:position}, revealed certain sensitivity of the optimization performance to an extra coefficient $\varepsilon$, referred to as \textit{timestep parameter} in \cite{Kaltsa}, according~to
\begin{align}\label{eq:positiondelta}
x_{i,d}^t=x_{i,d}^{t-1}+v_{i,d}^{t}+\varepsilon a_{i,d}^t.
\end{align}
In~Fig.~\ref{fig:apsofurther}, the performance of A-APSO taking into account \eqref{eq:positiondelta} with $\varepsilon=2$ instead of \eqref{eq:position}, i.e., $\varepsilon=0.5$, is shown to be even better when compared to the standard PSO algorithm.

\balance

\section{Conclusion}\label{sec:concl}

In this paper, we presented an improved version of the standard PSO algorithm, called A-APSO, which considers also the acceleration factor of the swarm particles in the associated model equations. The novel swarm intelligence algorithm has been employed for dynamic spectrum sensing in a cooperative cognitive radio network. The simulation results revealed that the proposed A-APSO is a promising tool for robust cooperative spectrum sensing, providing higher probabilities of detection than the standard PSO algorithm.

\section*{Acknowledgment}

The authors are grateful to Prof.~Traianos~V.~Yioultsis, Aristotle University of Thessaloniki, Department of Electrical and Computer Engineering, for his useful discussions.


\begin{thebibliography}{99}

\bibitem{Suraweera}
H. A. Suraweera, P. J. Smith, and M. Shafi, ``Capacity limits and performance analysis of cognitive radio with imperfect channel knowledge," \emph{IEEE Trans. Veh. Technol.}, vol.~59, no.~4, pp.~1811--1822, May 2010.

\bibitem{Nijsure}
Y. Nijsure, Y. Chen, S. Boussakta, C. Yuen, Y. H. Chew; Z. Ding, ``Novel system architecture and waveform design for cognitive radar radio networks," \emph{IEEE Trans. Veh. Technol.}, vol.~61, no.~8, pp.~3630--3642, Oct. 2012.

\bibitem{Stotas1}
S. Stotas and A. Nallanathan, ``On the throughput and spectrum sensing enhancement of opportunistic spectrum access cognitive radio networks," \emph{IEEE Trans. Wireless Commun.}, vol.~11, no.~1, pp.~97--107, Jan. 2012.

\bibitem{Stotas2}
S. Stotas and A. Nallanathan, ``Optimal sensing time and power allocation in multiband cognitive radio networks," \emph{IEEE Trans. Commun.}, vol.~59, no.~1, pp.~226--235, Jan. 2011.

\bibitem{Hanzo_2017}
W. Liang, S. X. Ng, and L. Hanzo, ``Cooperative overlay spectrum access in cognitive radio networks," \emph{IEEE Commun. Surveys Tuts.}, vol.~19, no.~3, pp.~1924--1944, 3rd Quart., 2017.

\bibitem{Chen_2014}
X. Chen, H.-H. Chen, and W. Meng, ``Cooperative communications for cognitive radio networks: From theory to applications," \emph{IEEE Commun. Surveys Tuts.}, vol.~16, no.~3, pp.~1180--1192, 3rd Quart., 2014.

\bibitem{daCosta}
D. B. da Costa, S. A\"{\i}ssa, and C. C. Cavalcante, ``Performance analysis of partial relay selection in cooperative spectrum sharing systems,'' \emph{Wirel. Pers. Commun.}, vol.~64, no.~1, pp.~79--92, May 2012.

\bibitem{Duong}
T. Q. Duong, D. B. da Costa, M. Elkashlan, and V. N. Q. Bao, ``Cognitive amplify-and-forward relay networks over Nakagami-$m$ fading," \emph{IEEE Trans. Veh. Technol.}, vol.~61, no.~5, pp.~2368--2374, Jun. 2012.

\bibitem{Zou}
Y. Zou, Y.-D. Yao, and B. Zheng, ``A cooperative sensing based cognitive relay transmission scheme without a dedicated sensing relay channel in cognitive radio networks,'' \emph{IEEE Trans. Signal Process.}, vol.~59, no.~2, pp.~854--858, Feb. 2011.

\bibitem{Lee}
C.-H. Lee and W. Wolf, ``Energy efficient techniques for cooperative spectrum sensing in cognitive radios,'' in \emph{Proc. IEEE Consumer Communications and Networking Conference}, Las Vegas, NV, Jan. 2008, pp.~968--972.

\bibitem{Quan}
Z. Quan, S. Cui, and A. H. Sayed, ``Optimal linear cooperation for spectrum sensing in cognitive radio networks," \emph{IEEE J. Sel. Topics Signal Process.}, vol.~2, no.~1, pp.~28--40, Feb. 2008.

\bibitem{Manesh_2017}
M. R. Manesh, A. Quadri, S. Subramaniam, and N. Kaabouch, ``An optimized SNR estimation technique using particle swarm optimization algorithm,'' in \emph{Proc. IEEE Computing and Communication Workshop and Conference}, Las Vegas, NV, Jan. 2017, pp.~1--6.

\bibitem{Alhammadi_2016}
A. Alhammadi, M. Roslee, M. Y. Alias, ``Analysis of spectrum handoff schemes in cognitive radio network using particle swarm optimization,'' in \emph{Proc. IEEE International Symposium on Telecommunication Technologies}, Kuala Lumpur, Malaysia, Nov. 2016, pp.~103--107.

\bibitem{Dawoud1_2016}
A. E. M. Dawoud, A. A. Rosas, M. Shokair, M. Elkordy, and S. El Halafawy, ``PSO-adaptive power allocation for multiuser GFDM-based cognitive radio networks,'' in \emph{Proc. International Conference on Selected Topics in Mobile \& Wireless Networking}, Cairo, Egypt, Apr. 2016, pp.~1--8.

\bibitem{Dawoud2_2016}
A. E. M. Dawoud, M. Shokair, M. Elkordy, and S. El Halafawy, ``Minimal throughput maximization for MIMO cognitive radio networks using particle swarm optimization,'' in \emph{Proc. International Conference on Selected Topics in Mobile \& Wireless Networking}, Cairo, Egypt, Apr. 2016, pp.~1--7.

\bibitem{Rashid}
R. A. Rashid, A. H. F. Bin Abdul Hamid, N. Fisal, S. K. Syed-Yusof, H. Hosseini, A. Lo, and A. Farzamnia, ``Efficient in-band spectrum sensing using swarm intelligence for cognitive radio network," \emph{Can. J. Elect. Comput. Eng.}, vol.~38, no.~2, pp.~106--115, Spr. 2015.

\bibitem{Pradhan}
P. M. Pradhan, G. Panda, and B. Majhi, ``Multiobjective cooperative spectrum sensing in cognitive radio using cat swarm optimization,'' in \emph{Proc. Wireless Advanced}, London, UK, Jun. 2012, pp.~44--48.

\bibitem{Zheng}
S. Zheng, C. Lou, and X. Yang, ``Cooperative spectrum sensing using particle swarm optimisation,'' \emph{Electron. Lett.}, vol.~46, no.~22, pp.~1525--1526, Oct. 2010.

\bibitem{Udgata}
S. K. Udgata, K. P. Kumar, and S. L. Sabat, ``Swarm intelligence based resource allocation algorithm for cognitive radio network," in \emph{Proc. International Conference on Parallel Distributed and Grid Computing}, Solan, India, Oct. 2010, pp.~324--329.

\bibitem{Kennedy}
J. Kennedy and R. Eberhart, ``Particle swarm optimization,'' in \emph{Proc. IEEE International Conference on Neural Networks}, Perth, Australia, Nov. 1995, pp.~1942--1948.

\bibitem{Shi1}
Y. Shi and R. C. Eberhart, ``Parameter selection in particle swarm optimization,'' in \emph{Proc. Conference on Evolutionary Programming}, San Diego, CA, Mar. 1998, pp.~591--600.

\bibitem{Shi2}
Y. Shi and R. C. Eberhart, ``A modified particle swarm optimizer,'' in \emph{Proc. IEEE International Conference on Evolutionary Computation}, Anchorage, AK, May 1998, pp.~69--73.

\bibitem{Han_2017}
H. Han, W. Lu, L. Zhang, and J. Qiao, ``Adaptive gradient multiobjective particle swarm optimization," \emph{IEEE Trans. Cybern.}, vol.~PP, no.~99, pp.~1--13, 2017.

\bibitem{Yang}
X.-S. Yang, S. Deb, and S. Fong, ``Accelerated particle swarm optimization and support vector machine for business optimization and applications,'' in \emph{Proc. International Conference on Networked Digital Technologies}, Macau, China, Jul. 2011, pp.~53--66.

\bibitem{Mohamed}
A. Z. Mohamed, S. H. Lee, H. Y. Hsu, and N. Nath, ``A faster path planner using accelerated particle swarm optimization,'' \emph{Artif. Life Robot.}, vol.~17, no.~2, pp.~233--240, Sep. 2012.

\bibitem{Elkaim}
G. H. Elkaim and M. Siegel, ``A lightweight control methodology for formation control of vehicle swarms,'' in \emph{Proc. World Congress of the International Federation of Automatic Control}, Prague, Czech Republic, Jul. 2005.

\bibitem{Apostolidis}
G. Apostolidis, ``Novel multiresolution signal analysis using swarm intelligence,'' Diploma Thesis, Aristotle University of Thessaloniki, Thessaloniki, Greece, Oct. 2012.

\bibitem{Kaltsa}
V. Kaltsa, A. Briassouli, I. Kompartsiaris, L. J. Hadjileontiadis, and M. G. Strintzis, ``Swarm intelligence for detecting interesting events in crowded environments,'' \emph{IEEE Trans. Image Process.}, vol.~24, no.~7, pp.~2153--2166, Jul. 2015.

\end{thebibliography}
\end{document}